# A HILBERT-MUMFORD CRITERION FOR $SL_2$-ACTIONS

JÜRGEN HAUSEN

ABSTRACT. Let the special linear group $G := SL_2$ act regularly on a $\mathbb{Q}$-factorial variety $X$. Consider a maximal torus $T \subset G$ and its normalizer $N \subset G$. We prove: If $U \subset X$ is a maximal open $N$-invariant subset admitting a good quotient $U \to U /\!\!/ N$ with a divisorial quotient space, then the intersection $W(U)$ of all translates $g \cdot U$ is open in $X$ and admits a good quotient $W(U) \to W(U) /\!\!/ G$ with a divisorial quotient space. Conversely, we obtain that every maximal open $G$-invariant subset $W \subset X$ admitting a good quotient $W \to W /\!\!/ G$ with a divisorial quotient space is of the form $W = W(U)$ for some maximal open $N$-invariant $U$ as above.

## INTRODUCTION

Given an action of a reductive group $G$ on a variety $X$, the task of Geometric Invariant Theory is the construction of open $G$-invariant subsets $W \subset X$ that admit reasonable quotients. We ask here for *good quotients*, that means $G$-invariant affine regular maps $p \colon W \to W /\!\!/ G$ of prevarieties such that the structure sheaf of $W /\!\!/ G$ equals the sheaf of invariants $p_*(\mathcal{O}_W)^G$, compare [9]. Note that we allow here nonseparated quotient spaces.

The task of Hilbert-Mumford Criteria is to reduce the construction of quotients $W \to W /\!\!/ G$ to the construction of quotients $U \to U /\!\!/ T$ for a maximal torus $T \subset G$. More precisely, one considers the following problem, compare [1], [2], [3] and [4]: Suppose that the open set $U \subset X$ is invariant under the normalizer $N \subset G$ of $T$ and admits a good quotient $U \to U /\!\!/ T$. When is the intersection $W(U)$ of all translates $g \cdot U$, $g \in G$, open in $X$ and admits a good quotient $W(U) \to W(U) /\!\!/ G$?

In this note we continue the study of the case $G = SL_2$ started in [2] and [3], where the above problem was solved for complete and for quasiprojective quotient spaces. Generalizing the latter setting, we focus on *divisorial* quotient spaces, i.e., prevarieties $Y$ such that each $y \in Y$ has an affine neighbourhood of the form $Y \setminus \mathrm{Supp}(D)$ with an effective Cartier divisor $D$ on $Y$, compare [6].

We shall use the following notions of maximality for open subsets with good quotient: Let $H \subset G$ be any reductive subgroup. By an *$H$-subset* we mean $H$-invariant open subset $U \subset X$ admitting a good quotient $U \to U /\!\!/ H$. We say that an $H$-subset $U \subset X$ is

- *d-maximal* if $U /\!\!/ H$ is divisorial, and $U$ does not occur as a saturated subset of a properly larger $H$-subset $U' \subset X$ with $U' /\!\!/ H$ divisorial,
- *s-d-maximal* if $U /\!\!/ H$ is separated, divisorial, and $U$ is not a saturated subset of a properly larger $H$-subset $U' \subset X$ with $U' /\!\!/ H$ separated and divisorial.

---

As before, fix a maximal torus $T \subset G$, and let $N \subset G$ be its normalizer. The aim of this note is to provide a recipe for constructing all d-maximal and all s-d-maximal $G$-subsets of a $G$-variety $X$ from the collection of all d-maximal $N$-subsets of $X$. The main result is the following, see Theorem 2.2:

**Theorem.** *Let $G = \mathrm{SL}_2(\mathbb{K})$ act regularly on a $\mathbb{Q}$-factorial variety $X$. For any open subset $U \subset X$ let $W(U)$ be the intersection of all translates $g \cdot U$, $g \in G$.*

(i) *If $U \subset X$ is a d-maximal $N$-subset then $W(U)$ is open in $X$ and admits a good quotient $W(U) \to W(U)/\!\!/G$ with a divisorial prevariety $W(U)/\!\!/G$.*
(ii) *Every d-maximal $G$-subset $W \subset X$ is of the form $W = W(U)$ for some d-maximal open $N$-subset $U \subset X$.*

Let us turn to the s-d-maximal $G$-subsets. By an *$N$-separated component* of a d-maximal $N$-subset $U' \subset X$ we mean the inverse image $U \subset U'$ of a maximal separated open subset of $U'/\!\!/N$ under the quotient map $U' \to U'/\!\!/N$. To any such $U \subset U'$, we associate its *$G$-kernel* $W(U)^\diamond$, see 3.1. This $G$-kernel admits a good quotient $W(U)^\diamond \to W(U)^\diamond /\!\!/ G$ with a separated divisorial quotient space. As a consequence of the main result, we obtain:

**Corollary.** *Every s-d-maximal $G$-subset $W \subset X$ is of the form $W = W(U)^\diamond$ with an $N$-separated component $U \subset U'$ of a d-maximal $N$-subset $U' \subset X$.*

Note that we always obtain algebraic (pre-)varieties as quotient spaces, whereas in [2] and [3] sincere algebraic spaces can occur. So, our results may also serve as algebraicity criteria for quotient spaces. The proof of the main result relies on the generalization of Mumford's Geometric Invariant Theory developped in [7]; see Section 1 for a summary and some slight extensions. The main result is proven in Section 2, and the proof of the corollary is performed in Section 3.

## 1. Generalized semistability

In [7] we generalized Mumford's construction of good quotients [8] by replacing his $G$-linearized line bundle with a certain group of Cartier divisors. The result is a theory producing all divisorial quotient spaces instead of only the quasiprojective ones. We recall here the basic results and adapt them to our actual purposes.

Let $X$ be an irreducible prevariety over an algebraically closed field $\mathbb{K}$ of characteristic zero. Fix a finitely generated free subgroup $\Lambda \subset \mathrm{CDiv}(X)$ of the group of Cartier divisors on $X$. Then one has an associated $\Lambda$-graded $\mathcal{O}_X$-algebra

$$\mathcal{A} := \bigoplus_{D \in \Lambda} \mathcal{A}_D := \bigoplus_{D \in \Lambda} \mathcal{O}_X(D).$$

The $\mathcal{O}_X$-algebra $\mathcal{A}$ gives rise to a prevariety $\widehat{X} := \mathrm{Spec}(\mathcal{A})$ and a canonical affine map $q \colon \widehat{X} \to X$ such that $q_*(\mathcal{O}_{\widehat{X}})$ equals $\mathcal{A}$. For a homogeneous local section $f \in \mathcal{A}_D(U)$, one defines its set of zeroes to be

$$Z(f) := \mathrm{Supp}(D|_U + \mathrm{div}(f)).$$

We call the group $\Lambda$ *ample* on an open subset $U \subset X$, if every $x \in U$ has an affine neighbourhood of the form $U \setminus Z(f)$ for a homogeneous $f \in \mathcal{A}(U)$. If $\Lambda$ is ample on $U$, then $\widehat{U} := q^{-1}(U)$ is a quasiaffine variety. Note that there exists a $\Lambda \subset \mathrm{CDiv}(X)$ which is ample on $X$ if and only if $X$ is *divisorial* in the sense of [6].

Now, let $G \times X \to X$ be a regular action of a reductive group on $X$. As in [7], we mean by a *G-linearization* of the group $\Lambda$ a graded $G$-sheaf structure on the $\mathcal{O}_X$-algebra $\mathcal{A}$ such that for every $G$-invariant open $U \subset X$ the induced representation of $G$ on $\mathcal{A}(U)$ is rational.

**Definition 1.1.** Let $\Lambda \subset \mathrm{CDiv}(X)$ be $G$-linearized, and let $U \subset X$ be a $G$-invariant open subset. A point $x \in U$ is called *$U$-semistable*, if $x$ has an affine neighbourhood $U' = U \setminus Z(f)$ with some $G$-invariant $f \in \mathcal{A}_D(U)$ such that the $D' \in \Lambda$ admitting a $G$-invariant $f' \in \mathcal{A}_{D'}(U')$ which is invertible in $\mathcal{A}(U')$ form a subgroup of finite index in $\Lambda$.

For a given $G$-linearized $\Lambda \subset \mathrm{CDiv}(X)$ and a $G$-invariant open $U \subset X$, the associated set of $U$-semistable points will be denoted by $U^{ss}(\Lambda)$, or $U^{ss}(\Lambda, G)$ if we want to specify the group $G$. Note that for $U = X$, the above definition specializes to the notion of semistability introduced in [7].

Following the lines of [7, Section 2], we show now that every set of $U$-semistable points admits a good quotient. First let us recall the precise definition:

**Definition 1.2.** A *good quotient* for the action of $G$ on $X$ is a $G$-invariant affine regular map $p \colon X \to X /\!\!/ G$ of prevarieties such that the canonical map $\mathcal{O}_{X /\!\!/ G} \to p_*(\mathcal{O}_X)^G$ is an isomorphism. A *geometric quotient* is a good quotient that separates orbits. Geometric quotients are denoted by $p \colon X \to X/G$.

Suppose now that $\Lambda \subset \mathrm{CDiv}(X)$ is $G$-linearized, and let $U \subset X$ be any $G$-invariant open subset. We shall need a geometric interpretation of $U$-semistability similar to the case $U = X$ treated in [7, Proposition 2.3]. For this we have to assume that $\Lambda$ is ample on $U$. Let $\mathcal{A}$ and $\widehat{X} := \mathrm{Spec}(\mathcal{A})$ be as before, and consider the canonical map $q \colon \widehat{X} \to X$.

Recall from [7, Section 1] that $q \colon \widehat{X} \to X$ is a geometric quotient for the action of the algebraic torus $H := \mathrm{Spec}(\mathbb{K}[\Lambda])$ on $\widehat{X}$ defined by the $\Lambda$-grading of $\mathcal{A}$. Moreover, the $G$-representation on $\mathcal{A}(U)$ induces a regular $G$-action on the quasiaffine variety $\widehat{U} = q^{-1}(U)$ such that the actions of $H$ and $G$ commute and $q \colon \widehat{U} \to U$ becomes $G$-equivariant.

Let $f_1, \ldots, f_r \in \mathcal{A}(U)$ be homogeneous and $G$-invariant such that the sets $U \setminus Z(f_i)$ are as in Definition 1.1 and cover $U^{ss}(\Lambda)$. Choose a $(G \times H)$-equivariant affine closure $\overline{U}$ of $\widehat{U}$ such that the functions $f_i \in \mathcal{O}(\widehat{U})$ extend regularly to $\overline{U}$ and $\overline{U}_{f_i} = \widehat{U}_{f_i}$ holds for each $i = 1, \ldots, r$. Then we have a good quotient
$$\overline{p} \colon \overline{U} \to \overline{U} /\!\!/ G := \mathrm{Spec}(\mathcal{O}(\overline{U}))^G.$$

The quotient variety $\overline{U} /\!\!/ G$ inherits a regular action of the torus $H$ such that the map $\overline{p} \colon \overline{U} \to \overline{U} /\!\!/ G$ becomes $H$-equivariant. Similar to [7, Proposition 2.3], we obtain:

**Lemma 1.3.** *Let $V_0 := \overline{U} /\!\!/ G \setminus \overline{p}(\overline{U} \setminus \widehat{U})$, and let $V_1 \subset \overline{U} /\!\!/ G$ be the union of all $H$-orbits with finite isotropy. Then*
$$q^{-1}(U^{ss}(\Lambda)) = \overline{p}^{-1}(V_0 \cap V_1).$$

*Proof.* Set $W := U^{ss}(\Lambda)$ and $\widehat{W} := q^{-1}(W)$. We begin with the inclusion "$\subset$". By [7, Remark 1.6] the set $q^{-1}(U \setminus Z(f_i))$ equals $\overline{U}_{f_i}$. Since each of the latter sets is $\overline{p}$-saturated and $\widehat{W}$ is covered by the $\overline{U}_{f_i}$, we see that $\widehat{W}$ is $\overline{p}$-saturated. In particular, we obtain $\overline{p}(\widehat{W}) \subset V_0$.

To verify $\overline{p}(\widehat{W}) \subset V_1$, let $z \in \widehat{W}$. Take one of the $f_i$ with $z \in \overline{U}_{f_i}$. As it is $G$-invariant, $f_i$ descends to an $H$-homogeneous function $h \in \mathcal{O}(\overline{U}/\!\!/G)$. By the properties of $f_i$, the function $h$ satisfies the condition of [7, Lemma 2.4] for the point $\overline{p}(z)$. Hence $H_{\overline{p}(z)}$ is finite, which means $\overline{p}(z) \in V_1$.

We come to the inclusion "⊃" of the assertion. Let $y \in V_0 \cap V_1$. Then [7, Lemma 2.4] provides an $h \in \mathcal{O}(\overline{U}/\!\!/G)$, homogeneous with respect to some $\chi^D \in \mathrm{Char}(H)$, such that $y \in V := (\overline{U}/\!\!/G)_h$ holds and the $D' \in \Lambda$ admitting an invertible $\chi^{D'}$-homogeneous function on $V$ form a subgroup of finite index in $\Lambda$. Suitably modifying $h$, we achieve additionally $V \subset V_0 \cap V_1$.

Now, consider a point $z \in \overline{p}^{-1}(y)$. Since $y \in V_0$, we have $z \in \widehat{U}$. We have to show that $q(z)$ is $U$-semistable. For this, consider the $G$-invariant homogeneous section $f := \overline{p}^*(h)|_{\widehat{U}}$ of $\mathcal{A}_D(U)$. By the choice of $h$, this $f$ fulfills the conditions of Definition 1.1 and thus the point $q(z)$ is in fact $U$-semistable. $\square$

As a consequence of this geometric description, we obtain existence of a good quotient for the set of $U$-semistable points for any $G$-invariant open subset $U \subset X$. The result generalizes [7, Theorem 3.1]:

**Proposition 1.4.** *Let $\Lambda \subset \mathrm{CDiv}(X)$ be $G$-linearized, and let $U \subset X$ be an open $G$-invariant subset. Then there is a good quotient*

$$p \colon U^{ss}(\Lambda) \to U^{ss}(\Lambda)/\!\!/G$$

*and the quotient space $U^{ss}(\Lambda)/\!\!/G$ is a divisorial prevariety. Moreover, for every $G$-invariant homogeneous $f \in \mathcal{A}(U^{ss}(\Lambda))$, the zero set $Z(f)$ is $p$-saturated.*

*Proof.* We may assume that $U = U^{ss}(\Lambda)$ holds. Then $\Lambda$ is ample on $U$, and we are in the setting of Lemma 1.3. Since the set $\widehat{U}$ is saturated with respect to the good quotient $\overline{p} \colon \overline{U} \to \overline{U}/\!\!/G$, restricting $\overline{p}$ to $\widehat{U}$ yields a good quotient $\widehat{p} \colon \widehat{U} \to \widehat{U}/\!\!/G$.

Moreover, Lemma 1.3 tells us that $H$ acts with at most finite isotropy groups on $\widehat{U}/\!\!/G$. Thus, there is a geometric quotient $\widehat{U}/\!\!/G \to (\widehat{U}/\!\!/G)/H$. By [7, Lemma 3.3], the latter quotient space is a divisorial prevariety. Since good quotients are categorical, we obtain a commutative diagram

$$\begin{array}{ccc} \widehat{U} & \xrightarrow{\widehat{p}} & \widehat{U}/\!\!/G \\ {\scriptstyle /H}\downarrow & & \downarrow{\scriptstyle /H} \\ U & \longrightarrow & (\widehat{U}/\!\!/G)/H \end{array}$$

Now it is straightforward to check that the induced map $U \to (\widehat{U}/\!\!/G)/H$ is the desired good quotient for the action of $G$ on $U$. This proves the first part of the assertion.

For the supplement, consider a $G$-invariant homogeneous $f \in \mathcal{A}(U^{ss}(\Lambda))$. By [7, Remark 1.6], the set $A$ of zeroes of $f \in \widehat{U}$ equals $q^{-1}(Z(f))$. Now, $A$ is $\widehat{p}$-saturated, and $\widehat{p}(A)$ is saturated with respect to the geometric quotient $\widehat{U}/\!\!/G \to (\widehat{U}/\!\!/G)/H$. Thus surjectivity of the involved maps gives the claim. $\square$

Similar to [8], we have also a converse of Proposition 1.4. Recall from [7, Section 4] that a group $\Lambda \subset \mathrm{CDiv}(X)$ is said to be *canonically $G$-linearized*, if on every homogeneous component $\mathcal{A}_D$ the $G$-sheaf structure arises from the action

$$(g \cdot f)(x) := f(g^{-1} \cdot x)$$

on the function field $\mathbb{K}(X)$. Now, assume that $X$ is a $\mathbb{Q}$-factorial $G$-variety, and let $U \subset X$ be an open $G$-invariant subset with a good quotient $U \to U/\!\!/G$. The proof of [7, Theorem 4.1] gives:

**Theorem 1.5.** *If $U/\!\!/G$ is divisorial, then there exists a canonically $G$-linearized $\Lambda \subset \mathrm{CDiv}(X)$ such that $U$ is contained in $X^{ss}(\Lambda)$ and is saturated with respect to the quotient map $X^{ss}(\Lambda) \to X^{ss}(\Lambda)/\!\!/G$.*

## 2. Proof of the main result

In this section, we prove our main result. First we have to introduce the following notions of maximality for open sets with good quotient, compare also [1]:

**Definition 2.1.** Let the reductive group $H$ act regularly on a variety $Y$. We say that $V \subset Y$ is an *$H$-subset* if it is open, $H$-invariant and admits a good quotient $V \to V/\!\!/H$. We say that an $H$-subset $V \subset Y$ is

   (i) *d-maximal* if $V/\!\!/H$ is divisorial, and $V$ does not occur as a saturated subset of a properly larger $H$-subset $V' \subset X$ with $V'/\!\!/H$ divisorial,
   (ii) *s-d-maximal* if $V/\!\!/H$ is separated and divisorial, and $V$ is not a saturated subset of a properly larger $H$-subset $V' \subset X$ with $V'/\!\!/H$ separated and divisorial.

Here a *saturated* subset of an $H$-subset $V'$ is a subset that is saturated with respect to the quotient map $V' \to V'/\!\!/H$.

Now, consider the special linear group $G := \mathrm{SL}_2(\mathbb{K})$. Fix a maximal torus $T \subset G$, and denote by $N$ its normalizer in $G$. For example:

$$T := \left\{ \begin{pmatrix} t & 0 \\ 0 & t^{-1} \end{pmatrix} ; t \in \mathbb{K}^* \right\}, \qquad N = T \cup nT, \quad \text{with } n := \begin{pmatrix} 0 & -1 \\ 1 & 0 \end{pmatrix}.$$

The main result of this note is generalizes and enhances [3, Theorem 9]. In the above notation, it reads as:

**Theorem 2.2.** *Let $G = \mathrm{SL}_2(\mathbb{K})$ act regularly on a $\mathbb{Q}$-factorial variety $X$. For any open subset $U \subset X$ let $W(U)$ be the intersection of all translates $g \cdot U$, $g \in G$.*

   (i) *If $U \subset X$ is a d-maximal $N$-subset then $W(U)$ is open and saturated in $U$, and there is a good quotient $W(U) \to W(U)/\!\!/G$ with a divisorial prevariety $W(U)/\!\!/G$.*
   (ii) *Every d-maximal $G$-subset $W \subset X$ is of the form $W = W(U)$ with a d-maximal $N$-subset $U \subset X$.*

In the proof we shall use the techniques presented in Section 1. Let $\Lambda \subset \mathrm{CDiv}(X)$ be a $G$-linearized group. Then $\Lambda$ is also linearized with respect to every subgroup of $G$. In particular, we obtain a set $U^{ss}(\Lambda, H)$ of semistable points for $G$-invariant open $U \subset X$ and every reductive subgroup $H \subset G$.

**Remark 2.3.** Let $U \subset X$ be a $G$-invariant open subset, let $H \subset G$ be a reductive subgroup, and let $g \in G$. Then a point $x \in U$ is $U$-semistable with respect to $H$ if and only if $g \cdot x$ is $U$-semistable with respect to $gHg^{-1}$.

The crucial step in the proof of Theorem 2.2 is to express generalized semistability in terms of maximal tori. For this, let $\mathrm{MT}(G)$ denote the set of maximal tori of the group $G$, and fix an element $T \in \mathrm{MT}(G)$.

**Lemma 2.4.** *Let $G = \mathrm{SL}_2(\mathbb{K})$ act regularly on a variety $X$, let $\Lambda \subset \mathrm{CDiv}(X)$ be a $G$-linearized group, and let $U \subset X$ be any $G$-invariant subset such that $\Lambda$ is ample on $U$. Then we have:*

$$U^{ss}(\Lambda, G) \;=\; \bigcap_{S \in \mathrm{MT}(G)} U^{ss}(\Lambda, S) \;=\; \bigcap_{g \in G} g \cdot U^{ss}(\Lambda, T).$$

*Moreover, the set $U^{ss}(\Lambda, G)$ is saturated in $U^{ss}(\Lambda, T)$ with respect to the quotient map $U^{ss}(\Lambda, T) \to U^{ss}(\Lambda, T) /\!/ T$.*

*Proof.* The supplement is due to Proposition 1.4. The second equality is clear by Remark 2.3. Moreover, the inclusion "$\subset$" of the first equality holds by the definition of semistability. Thus we are left with proving the inclusion "$\supset$" of the first equality.

For this, let $\mathcal{A}$ be the graded $\mathcal{O}_X$-algebra associated to $\Lambda$, and set $\widehat{X} := \mathrm{Spec}(\mathcal{A})$. Moreover, let $q \colon \widehat{X} \to X$ be the canonical map and denote by $H := \mathrm{Spec}(\mathbb{K}[\Lambda])$ the algebraic torus acting on $\widehat{X}$. Finally, let $\widehat{U} := q^{-1}(U)$.

Choose $G$-invariant homogeneous $f_1, \ldots, f_r \in \mathcal{A}(U)$ and $T$-invariant homogeneous $h_1, \ldots, h_s \in \mathcal{A}(U)$ such that the complements $U \setminus Z(f_i)$ and $U \setminus Z(h_j)$ satisfy the conditions of Definition 1.1 and

$$\begin{aligned} U^{ss}(\Lambda, G) &= (U \setminus Z(f_1)) \cup \ldots \cup (U \setminus Z(f_r)), \\ U^{ss}(\Lambda, T) &= (U \setminus Z(h_1)) \cup \ldots \cup (U \setminus Z(h_r)). \end{aligned}$$

Since $\Lambda$ is ample on $U$, there is a $(G \times H)$-equivariant affine closure $\overline{U}$ of $\widehat{U}$ such that the $f_i$ and $h_j$ extend to regular functions on $\overline{U}$ satisfying $\overline{U}_{f_i} = \widehat{U}_{f_i}$ and $\overline{U}_{h_j} = \widehat{U}_{h_j}$. Moreover, we obtain a commutative diagram of $H$-equivariant maps:

$$\begin{array}{ccc} \overline{U} & \xrightarrow{\overline{p}_G} & \overline{U}/\!/G \\ {\scriptstyle \overline{p}_T} \searrow & {\scriptstyle /\!/T} \;\; {\scriptstyle /\!/G} \nearrow & \\ & \overline{U}/\!/T & \end{array}$$

In the further proof, we shall apply the geometric characterization of semistability given in Lemma 1.3 to this diagram. A first step is to verify

$$\bigcap_{S \in \mathrm{MT}(G)} q^{-1}(U^{ss}(\Lambda, S)) \;\subset\; \widehat{U} \setminus \overline{p}_G^{-1}(\overline{p}_G(\overline{U} \setminus \widehat{U})).$$

Let $x \in U^{ss}(\Lambda, S)$ for all maximal tori $S \subset G$, and let $z \in q^{-1}(x)$. Set $y := \overline{p}_G(z)$. Suppose $y \in \overline{p}_G(\overline{U} \setminus \widehat{U})$. Let $Gz'$ be the closed orbit in $\overline{p}_G^{-1}(y)$. Then $Gz'$ is contained in $\overline{U} \setminus \widehat{U}$. Moreover, the Hilbert-Mumford-Birkes Lemma [5], provides a maximal torus $S \subset G$ such that the closure of $S \cdot z$ intersects $G \cdot z'$.

Let $g \in G$ with $gSg^{-1} = T$. Then the closure of $T \cdot g \cdot z$ contains a point $z'' \in G \cdot z'$. Surely, $\overline{p}_T(g \cdot z)$ equals $\overline{p}_T(z'')$. Thus, since $z'' \in \overline{U} \setminus \widehat{U}$, Lemma 1.3 tells us that $g \cdot x = q(g \cdot z)$ is not semistable with respect to $T$. By Remark 2.3, the point $x$ is not semistable with respect to $S$. A contradiction. So the first step is done.

Now assume that there is a point $x \in X$ that is semistable with respect to all maximal tori of $G$ but not with respect to $G$ itself. By the above inclusion, we have $x = q(z)$ for some $z \in \widehat{U} \setminus \overline{p}_G^{-1}(\overline{p}_G(\overline{U} \setminus \widehat{U}))$. Lemma 1.3 tells us that for $y := \overline{p}_G(z)$ the isotropy group $H_y$ is infinite. Let $H_0$ be the connected component of the neutral element of $H_y$, and let $G \cdot z_0$ be the closed $G$-orbit in the fibre $\overline{p}_G^{-1}(y)$.

Then $H_0$ acts freely on the fibre $\overline{p}_G^{-1}(y) \subset \widehat{U}$. Since $G \cdot z_0$ is the only closed $G$-orbit in this fibre, it is invariant under the action of $H_0$. Let $G_0$ denote the stabilizer of $H_0 \cdot z_0$. We claim that $G_0 \cdot z_0 = H_0 \cdot z_0$ holds. Indeed, let $h \in H_0$. Since $H_0 \cdot z_0 \subset G \cdot z_0$, there is a $g \in G$ with $g \cdot z_0 = h \cdot z_0$. We have $g \in G_0$ because

$$g \cdot (H_0 \cdot z_0) = H_0 \cdot (g \cdot z_0) = H_0 \cdot (h \cdot z_0) = H_0 \cdot z_0.$$

Denoting by $\mu \colon G_0 \to G_0 \cdot z_0$ and $\tau \colon H_0 \to H_0 \cdot z$ the orbit maps, we obtain an epimorphism $G_0 \to H_0$, $g \mapsto \tau^{-1}(\mu(g))$ of algebraic groups. In particular, for every maximal torus $S_0 \subset G_0$ we have $S_0 \cdot z_0 = H_0 \cdot z_0$. Consequently, the maximal tori of $G_0$ are nontrivial, and hence of dimension one. So, applying Lemma 1.3 to a maximal torus of $G_0$ yields

$$q(z_0) \notin \bigcap_{S \in \mathrm{MT}(G)} U^{ss}(\Lambda, S).$$

By the choice of $z$, this implies that $G \cdot z \neq G \cdot z_0$. Hence $G \cdot z_0$ is of smaller dimension than $G \cdot z$. In other words, the isotropy group $G_{z_0}$ is infinite. Note that $G_{z_0} \subset G_0$. Since $G \cdot z_0$ is affine, $G_{z_0}$ is reductive and hence contains a nontrivial torus $S_0$ of $G_0$. Since $S_0$ is by dimension reasons already a maximal torus of $G_0$, this contradicts $S_0 \cdot z_0 = H_0 \cdot z_0$. □

*Proof of Theorem 2.2.* We begin with (i). By Proposition 1.5, there is a canonically $N$-linearized group $\Lambda \subset \mathrm{CDiv}(X)$ such that $U = X^{ss}(\Lambda, N)$ holds. We show that on a subgroup $\Lambda'' \subset \Lambda$ of finite index, the canonical $N$-linearization of $\Lambda''$ extends to a $G$-linearization of $\Lambda''$, compare also [7, Proof of Theorem 5.1]:

Using Sumihiro's Equivariant Completion Theorem [10, Theorem 3] and equivariant normalization, we find a complete normal $G$-variety $Y$ which contains $X$ as a $G$-invariant open subset. Consider the set of singular points $Y_{\mathrm{sing}} \subset Y$, and remove the closed $G$-invariant set $Y_{\mathrm{sing}} \setminus X$ from $Y$. Then $Y$ is possibly no longer complete, but by normality, we still have $\mathcal{O}(Y) = \mathbb{K}$.

By construction, every Cartier divisor $D$ on $X$ extends to a Cartier divisor $E$ on $Y$: just replace the components of $D$ with their closures in $Y$. The resulting Weil divisor $E$ on $Y$ is Cartier, because its restriction $X$ is so, and any point of $Y \setminus X$ is smooth. Proceeding this way, we can extend the group $\Lambda \subset \mathrm{CDiv}(X)$ to a canonically $N$-linearized group $\Gamma \subset \mathrm{CDiv}(Y)$.

On the other hand, some subgroup $\Gamma' \subset \Gamma$ of finite index admits a $G$-linearization. Since we have $\mathcal{O}(Y) = \mathbb{K}$, we can apply [7, Proposition 1.5] to see that on some further subgroup $\Gamma'' \subset \Gamma$ of finite index the $N$-linearization inherited by the $G$-linearization and the canonical $N$-linearization coincide. Thus, restricting $\Gamma''$ to $X$ gives the desired subgroup $\Lambda'' \subset \Lambda$.

We replace $\Lambda$ with the above $\Lambda''$. Note that this does not affect $U = X^{ss}(\Lambda, N)$. Let $\mathcal{A}$ denote the graded $\mathcal{O}_X$-algebra defined by $\Lambda$, and let $V \subset X$ consist of all points admitting an affine neighbourhood $X \setminus Z(f)$ with a homogeneous section $f \in \mathcal{A}(X)$. Then $V$ is open and $G$-invariant, and the group $\Lambda$ is ample on $V$.

By the definition of semistability, we have $U \subset V^{ss}(\Lambda, N)$. Moreover, since $U$ is defined by removing zero sets of global $N$-invariant homogeneous sections, the supplement of Proposition 1.4 tells us that $U$ is even saturated with respect to the quotient map $V^{ss}(\Lambda, N) \to V^{ss}(\Lambda, N) /\!/ N$. Since $U$ is a d-maximal $N$-set, we have

$U = V^{ss}(\Lambda, N)$. So, Lemma 2.4 yields

$$W(U) = \bigcap_{g \in G} g \cdot U = \bigcap_{g \in G} g \cdot V^{ss}(\Lambda, N) \subset \bigcap_{g \in G} g \cdot V^{ss}(\Lambda, T) = V^{ss}(\Lambda, G).$$

Since $V^{ss}(\Lambda, G)$ is contained in $V^{ss}(\Lambda, N)$, we have in fact equality. In particular, $W(U)$ is open, and Proposition 1.4 yields a good quotient $W(U) \to W(U)/\!/G$ with a divisorial quotient space $W(U)/\!/G$. Moreover, we see that $W(U)$ is saturated with respect to the quotient map $U \to U/\!/N$.

We prove (ii). If $W \subset X$ is a d-maximal $G$-set, then Proposition 1.5 tells us that $W = X^{ss}(\Lambda, G)$ holds for some $G$-linearized group $\Lambda \subset \mathrm{CDiv}(X)$. Let $U \subset X$ be any d-maximal $N$-set $U \subset X$ containing $X^{ss}(\Lambda, N)$ as a saturated open subset. By assertion (i), the set $W(U)$ admits a good quotient by the action of $G$ with a divisorial quotient space $W(U)/\!/G$.

Surely, we have $W \subset W(U)$. Moreover, $W$ is saturated in $U$, and $W(U)$ is saturated in $U$. Thus we conclude that $W$ is saturated in $W(U)$ with respect to the quotient map $W(U) \to W(U)/\!/N$. The classical Hilbert-Mumford Lemma [5] yields that $W$ is even saturated with respect to $W(U) \to W(U)/\!/G$, see also [4, Proposition 2.6]. By d-maximality of $W$, this means $W = W(U)$. □

## 3. Proof of the Corollary

As in the previous section, $G = \mathrm{SL}_2$ acts on a $\mathbb{Q}$-factorial variety $X$, and $N \subset G$ denotes the normalizer of a maximal torus $T \subset G$. Let us give the precise definition of the $G$-kernel:

Let $U' \subset X$ be a d-maximal $N$-subset. Then we have good quotients $s \colon U' \to U'/\!/N$ and $p \colon W(U') \to W(U')/\!/G$, where the latter exists by Theorem 2.2. By an *$N$-separated component* $U \subset U'$ we mean the inverse image $U = s^{-1}(V)$ of a maximal separated open subset $V \subset U'/\!/N$.

**Definition 3.1.** The *$G$-kernel* $W(U)^\diamond$ of an $N$-separated component $U \subset U'$ is defined as follows: Let $W_0$ consist of all $x \in W(U)$ with $p^{-1}(p(x)) \subset W(U)$. Then $W(U)^\diamond$ is the inverse image $p^{-1}(V_0)$ of the set $V_0$ of interior points of $p(W_0)$.

Note that a $G$-kernel can be empty. As an application of Theorem 2.2 we show now that every s-d-maximal $G$-subset of $X$ is a $G$-kernel:

**Corollary 3.2.**   (i) *Let $U \subset U'$ be an $N$-separated component of a d-maximal $N$-subset $U' \subset X$. Then there is a good quotient $W(U)^\diamond \to W(U)^\diamond/\!/G$ with a separated divisorial quotient space.*
  (ii) *Every s-d-maximal $G$-subset $W \subset X$ is of the form $W = W(U)^\diamond$ with an $N$-separated component $U \subset U'$ of a d-maximal $N$-subset $U' \subset X$.*

*Proof.* In order to prove (i), note first that the $G$-kernel $W(U)^\diamond$ is by construction saturated with respect to the quotient map $W(U') \to W(U')/\!/G$. Hence $W(U)^\diamond$ is also saturated with respect to the quotient map $W(U') \to W(U')/\!/N$. Since $W(U') \subset U'$ and $U \subset U'$ are saturated inclusions of $N$-subsets, we obtain that $W(U)^\diamond$ is saturated in the $N$-subset $U$.

Consequently, $W(U)^\diamond$ admits a good quotient $W(U)^\diamond \to W(U)^\diamond/\!/N$ with a separated divisorial quotient space. Thus we may use for example [7, Theorem 5.1] to infer existence of a good quotient $W(U)^\diamond \to W(U)^\diamond/\!/G$ with a separated divisorial quotient space.

We prove (ii). Let $W \subset X$ be an s-d-maximal $G$-subset. Then $W$ is a saturated subset of some d-maximal $G$-subset $W' \subset X$. We consider the following commutative diagram:

$$\begin{array}{ccc} W' & \xrightarrow{p} & W'/\!\!/G \\ & \searrow{\scriptstyle q} \quad \nearrow{\scriptstyle r} & \\ & W'/\!\!/N & \end{array}$$

Note that $Z := p(W)$ is a maximal separated open subset of $W'/\!\!/G$. Since $r$ is an affine map, the inverse image $r^{-1}(Z)$ is contained in some maximal separated open subset $Y \subset W'/\!\!/N$.

By Theorem 2.2, we have $W' = W(U')$ with some d-maximal $N$-subset $U' \subset X$. Moreover, $W'$ is saturated with respect to the quotient map $s \colon U' \to U'/\!\!/N$. So, we may view $q$ as the restriction of $s$. Then $Y$ is of the form $s(U) \cap s(W')$ with some $N$-separated component $U \subset U'$. We have

$$W(U) = W(U) \cap W' = W(U \cap W') = W(q^{-1}(Y)).$$

Consequently, $W = p^{-1}(Z)$ is contained in $W(U)$, and hence in the $G$-kernel $W(U)^\diamond$. According to (i), there is a good quotient $W(U)^\diamond \to W(U)^\diamond/\!\!/G$ with a separated divisorial quotient space. Moreover, $W \subset W'$ and $W(U)^\diamond \subset W'$ are saturated inclusions of $G$-sets. Thus $W \subset W(U)^\diamond$ is saturated, and, by s-d-maximality, we have $W = W(U)^\diamond$. $\square$

## References

[1] A. Białynicki-Birula: Finiteness of the number of maximal open subsets with good quotients. Transform. Groups Vol. 3, No. 4, 301–319 (1998)

[2] A. Białynicki-Birula, J. Święcicka: On complete orbit spaces of SL(2)-actions. Colloq. Math. Vol. 55, No. 2, 229–2431 (1988)

[3] A. Białynicki-Birula, J. Święcicka: On complete orbit spaces of SL(2)-actions II. Colloq. Math. Vol. 63, No. 1, 9–20 (1992)

[4] A. Białynicki-Birula, J. Święcicka: Open subsets of projective spaces with a good quotient by an action of a reductive group. Transform. Groups 1, 153–185 (1996)

[5] D. Birkes: Orbits of linear algebraic groups. Ann. Math., II. Ser. 93, 459–475 (1971)

[6] M. Borelli: Divisorial varieties. Pacific J. Math. 13, 375–388 (1963)

[7] J. Hausen: A Generalization of Mumford's Geometric Invariant Theory. Documenta Math. 6, 571–592 (2001)

[8] D. Mumford, J. Fogarty, F. Kirwan: Geometric Invariant Theory; third enlarged edition. Springer, Berlin, Heidelberg, 1994

[9] C.S. Seshadri: Quotient spaces modulo reductive algebraic groups. Ann. of Math. 95, 511–556 (1972)

[10] H Sumihiro: Equivariant completion. J. Math. Kyoto Univ. 14, 1–28 (1974)

Fachbereich Math. und Statistik, Universität Konstanz, 78457 Konstanz, Germany
*E-mail address*: Juergen.Hausen@uni-konstanz.de